\documentstyle[12pt]{article}

\input amssym.def
\input amssym.tex

\def\vv{{\underline{v}}}
\def\tt{{\underline{t}}}
\def\aa{\underline{a}}

\def\mm{\underline{m}}
\def\kk{\underline{k}}
\def\1{\underline{1}}
\def\R{\Bbb R}
\def\P{\Bbb P}

\def\L{\Bbb L}
\def\Z{\Bbb Z}

\def\C{\Bbb C}
\def\CP{\Bbb C\Bbb P}

\newtheorem{theorem}{Theorem}
\newtheorem{statement}{Statement}
\newtheorem{corollary}{Corollary}
\newtheorem{lemma}{Lemma}

\newenvironment{definition}
{\smallskip\noindent{\bf Definition\/}:}{\smallskip\par}

\newenvironment{remark}
{\smallskip\noindent{\bf Remark\/}.}{\smallskip\par}
\newenvironment{remarks}
{\smallskip\noindent{\bf Remarks\/}.}{\smallskip\par}

\title{Integrals with respect to the Euler characteristic over
spaces of functions and the Alexander polynomial.}
\author{A.Campillo
\and F.Delgado \thanks{First two authors were partially supported by
DGICYT PB97-0471. Address:
University of Valladolid, Dept. of Algebra, Geometry and Topology,
47005 Valladolid, Spain.
E-mail: campillo\symbol{'100}cpd.uva.es, fdelgado\symbol{'100}agt.uva.es}
\and S.M.Gusein-Zade \thanks{Partially supported by the grants
RFBR, INTAS--97--1644 and NWO 047.008.005.
Address: Moscow State University,
Dept. of Mathematics and Mechanics, Moscow, 119899, Russia.
E-mail: sabir\symbol{'100}mccme.ru}}
\date{}
\begin{document}

\def\eps{\varepsilon}

\maketitle


Recently there were obtained some formulae which express the Alexander
polynomial (and thus the zeta-function of the classical monodromy
transformation) of a plane
curve singularity in terms of the ring of functions on the curve:
\cite{CDG4}. One of them describes the coefficients of the Alexander
polynomial or of zeta-function of the monodromy transformation
as Euler characteristics
of some explicitly constructed spaces. For the Alexander polynomial these
spaces are complements to arrangements of projective hyperplanes in
projective spaces. Thus for the zeta-function they are disjoint unions
of such spaces. A little bit later J.Denef and F.Loeser (\cite{DL})
described the
Lefschetz numbers of iterates of the monodromy transformation of a
hypersurface singularity (of any dimension) as Euler characteristics of
some subspaces of the space of (truncated) arcs. After that it was understood
that the results of \cite{CDG4} are connected with the notion of the motivic
integration or rather with its version (in some sense a dual one) where
the space of arcs is substituted by the space of functions. The final
aim of this paper is to discuss the notion of the integral with respect
to the Euler characteristics (or with respect to the generalized Euler
characteristic) over the space of functions (or over its projectivization)
and its connection with the formulae for the coefficients of the Alexander
polynomial and of the zeta-function of the monodromy transformation as Euler
characteristics of some spaces.

There were two results which preceded the formulae for the Alexander
polynomial of a plane curve singularity from \cite{CDG4}. At first there
was seen no real connection between them. We discuss these results in the
next two sections.

\section{The Poincar\'e series of the ring of functions and the zeta-function
of the monodromy transformation}\label{sec1}

One of the results (or rather an observation) identifies the Poincar\'e
series of
the ring of functions on an irreducible plane curve singularity with the
zeta-function of the classical monodromy transformation of it (\cite{CDG2}).

Let $(C, 0)\subset(\C^2, 0)$ be an irreducible germ of a plane curve.
Let $f=0$ be an equation of the curve $C$ ($f:(\C^2, 0)\to (\C, 0)$
is a germ of a holomorphic function with an isolated critical point
at the origin), let $\varphi:(\C, 0)\to(\C^2, 0)$ be a parameterization
(uniformization) of the curve $C$, i.e., a germ of an analytic map such
that ${\rm{Im}}\,\varphi=C$ and $\varphi$ is an isomorphism between $\C$
and $C$ outside of the origin.
For a germ $g$ from the ring ${\cal O}_{\C^2,0}$ of germs of holomorphic
functions at the origin in the plane $\C^2$ (or from the ring
${\cal O}_{C,0}={\cal O}_{C^2,0}/(f)$
of functions on the curve $C$), let $v(g)\in\Z_{\ge 0}$ be the order of
zero at the origin of the function $g\circ\varphi:(\C, 0)\to \C$, i.e.,
$g\circ\varphi(\tau)=a\cdot\tau^{v(g)}+{\ terms~of~higher~degree}$, $a\ne 0$
(if $g\circ\varphi\equiv 0$, i.e., if $g\in(f)$, one assumes $v(g)=+\infty$).
The set $S_C$ of integers $v=v(g)$ for all germs $g\in{\cal O}_{\C^2,0}$
with $v(g)<+\infty$ is a subsemigroup in $\Z_{\ge 0}$ and is called the
semigroup (of values) of the (irreducible) plane curve singularity $C$.

Let $\varphi^*: {\cal O}_{C,0}\to {\cal O}_{\C, 0}$ be the natural
homomorphism induced by the map $\varphi$. The ring ${\cal O}_{\C, 0}$ of
germs of holomorphic functions in one variable has a natural
filtration by powers of the maximal ideal ${\goth m}\subset{\cal O}_{\C, 0}$:
${\cal O}_{\C, 0}={\goth m}^0\supset{\goth m}^1\supset{\goth m}^2\supset\ldots$
Let $W_i=(\varphi^*)^{-1}({\goth m}^i)\subset {\cal O}_{C,0}$. The ideals
$W_i$ form a filtration of the ring ${\cal O}_{C,0}$ of functions on the curve
$C$: ${\cal O}_{C, 0}=W_0\supset W_1\supset W_2\supset\ldots$ Since
$\dim {\goth m}^i/{\goth m}^{i+1}=1$, the dimension of each factor $W_i/W_{i+1}$
is either $1$ or $0$. It is equal to $1$ if and only if $i$ is an element of
the semigroup $S_C$ of the curve $C$.

Let $P_C(t)=\sum\limits_{i=0}^\infty \dim(W_i/W_{i+1})\cdot t^i$
be the Poincar\'e series of the filtration $W_\bullet$.
One can write $P_C(t)$ as $\sum\limits_{v\in S_C} t^v$.
Since all sufficiently large integers belong to the semigroup $S_C$,
$P_C(t)$ is the power series decomposition of a rational function (which we
also denote by $P_C(t)$), moreover of a polynomial divided by $(1-t)$.

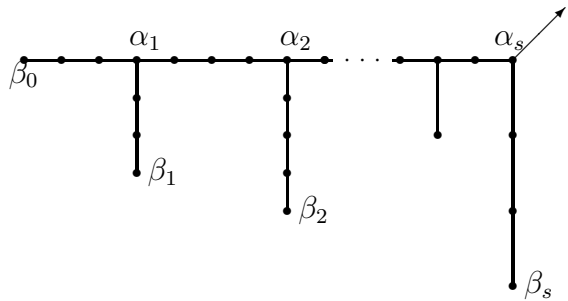
\begin{figure}
$$
\unitlength=1.00mm
\begin{picture}(80.00,30.00)(-10,3)
\thicklines
\put(-5,30){\line(1,0){41}}
\put(44,30){\line(1,0){16}}
\put(38,30){\circle*{0.5}}
\put(40,30){\circle*{0.5}}
\put(42,30){\circle*{0.5}}
\put(30,10){\line(0,1){20}}
\put(50,20){\line(0,1){10}}
\put(60,0){\line(0,1){30}}
\put(10,15){\line(0,1){15}}
\thinlines
\put(60,30){\vector(1,1){7}}
\put(20,30){\circle*{1}}
\put(30,30){\circle*{1}}
\put(50,30){\circle*{1}}
\put(60,30){\circle*{1}}
\put(30,20){\circle*{1}}
\put(60,20){\circle*{1}}
\put(60,10){\circle*{1}}
\put(10,30){\circle*{1}}
\put(30,10){\circle*{1}}
\put(50,20){\circle*{1}}
\put(60,0){\circle*{1}}
\put(-5,30){\circle*{1}}
\put(0,30){\circle*{1}}
\put(5,30){\circle*{1}}
\put(15,30){\circle*{1}}
\put(25,30){\circle*{1}}
\put(35,30){\circle*{1}}
\put(45,30){\circle*{1}}
\put(55,30){\circle*{1}}
\put(10,25){\circle*{1}}
\put(10,20){\circle*{1}}
\put(10,15){\circle*{1}}
\put(30,25){\circle*{1}}
\put(30,15){\circle*{1}}
\put(35,30){\circle*{1}}
\put(-7,27){$\beta_0$}
\put(11.5,14){$\beta_1$}
\put(31.5,9){$\beta_2$}
\put(61.5,-1){$\beta_s$}
\put(9,32){$\alpha_1$}
\put(29,32){$\alpha_2$}
\put(57.5,32){$\alpha_s$}
\end{picture}
$$
\caption{The dual resolution graph of an irreducible curve $C$.}
\label{fig1}
\end{figure}

Combinatorial properties of the semigroup of an irreducible plane curve
singularity (see, e.g., \cite{ZT}) permit to calculate the rational function
$P_C(t)$ in the following terms. Let $\pi:({\cal X}, {\cal D})\to (\C^2, 0)$
be the minimal embedded resolution of the plane curve singularity $C$. The
exceptional divisor ${\cal D}=\pi^{-1}(0)$ of the resolution is the union
of its irreducible components, each of them isomorphic to the projective
line $\C\P^1$. A resolution $\pi:({\cal X}, {\cal D})\to (\C^2, 0)$
of a plane curve singularity $C=\{f=0\}$ (not necessarily irreducible)
can be described by its dual
graph $\Gamma$. Vertices of the graph $\Gamma$ are in one--to--one
correspondence with irreducible components of the total transform
$\pi^{-1}(C)$ of the curve $C$, i.e., with components of the exceptional
divisor ${\cal D}$ and of the strict transform $\overline{\pi^{-1}(C)
\setminus{\cal D}}$ of the curve $C$. The vertices corresponding
to the components of the strict transform of the curve $C$ are depicted by
arrows. Two vertices of the graph $\Gamma$ are connected by an edge iff
the corresponding components of the total transform $\pi^{-1}(C)$
intersect. For a vertex $\sigma$ corresponding to a component
$E_\sigma$ of the exceptional divisor ${\cal D}$, let $m^\sigma$
be the multiplicity (order) of the lifting $f\circ\pi$ of the
function $f$ (the equation of the curve $C$) along $E_\sigma$.
The dual graph of the minimal resolution of an irreducible plane curve
singularity $C$ is shown in Fig.{\ref{fig1}}.
Here $s$ is the number of Puiseux pairs of the curve $C=\{f=0\}$. Let
$\beta_i$ ($i=0, 1, \ldots, s$) and $\alpha_i$ ($i=1, 2, \ldots, s$)
be the dead
ends and the star points of the graph $\Gamma$ (see Fig.{\ref{fig1}}), let
$\overline{\beta}_j=m^{\beta_j}$, $\overline{\alpha}_j=m^{\alpha_j}$.
The set of integers $\{\overline{\beta}_j\vert
j=0, 1, \ldots, s\}$ is the minimal set of generators of the semigroup $S_C$
of the curve $C$. It is known that the integers $\overline{\alpha}_j$
are multiples of the integers $\overline{\beta}_j$ for $j=1, \ldots, s$:
$\overline{\alpha}_j=(n_j+1)\overline{\beta}_j$, and that
each element $v\in S_C$ can be represented in a unique way in the form
$v=k_0\overline{\beta}_0+\sum\limits_{j=1}^s k_j\overline{\beta}_j$ with
$k_0\ge 0$, $0\le k_j\le n_j$ for $1\le j\le s$. Using this fact one can
obtain the following formula (see, e.g., \cite{CDG2}).

\begin{statement}
\begin{equation}\label{eq1}
P_C(t)=\frac
{\prod\limits_{j=1}^s(1-t^{\overline{\alpha}_j})}
{\prod\limits_{j=0}^s(1-t^{\overline{\beta}_j})}.
\end{equation}
\end{statement}

For a plane curve singularity $C$ (again not necessarily irreducible)
defined by an equation $\{f=0\}$ let $h_f:V_f\to V_f$ be the classical
monodromy transformation of the germ $f$. Here $V_f$ is the Milnor fibre of the
singularity $f$: $V_f=\{z\in\C^2: \Vert z\Vert\le\varepsilon, f(z)=\delta\}$,
where $0<\vert\delta\vert\ll\varepsilon$, $\varepsilon$ is small enough.
The monodromy transformation $h_f$ is a diffeomorphism of the Milnor fibre
$V_f$ well defined up to isotopy. The zeta-function $\zeta_h(t)$ of a
transformation $h:X\to X$ is the rational function in the variable $t$
defined as follows:
$$\zeta_h(t)=
\prod\limits_{q\geq 0} \left[\det\left(id-t\cdot {h_*}|_{H_q(X;\R)}
\right)\right]^{(-1)^{q+1}}.$$
Let $\zeta_C(t)=\zeta_{h_f}(t)$ be the zeta-function of the classical
monodromy transformation $h_f$. By the formula of N.A'Campo (\cite{A'C})
the zeta-function $\zeta_C(t)$ of a plane curve singularity $C=\{f=0\}$
can be expressed in terms of an embedded resolution $\pi:({\cal X},
{\cal D})\to (\C^2, 0)$ of the curve $C$ as
\begin{equation}\label{eq2}
\zeta_C(t)=\prod\limits_{E_\sigma\subset{\cal D}}
\left(1-t^{m^\sigma}\right)^{-\chi({\stackrel{\circ}{E}}\sigma)},
\end{equation}
where ${\stackrel{\circ}{E}}_\sigma$ is "the smooth part" of the component
$E_\sigma$ of the exceptional divisor ${\cal D}$, i.e., $E_\sigma$ minus
intersection points with all other components of the total transform of
the curve $C$.

For the irreducible plane curve singularity $C$ (the dual graph $\Gamma$
of the minimal resolution of which is shown in Fig.{\ref{fig1}})
$\chi({\stackrel{\circ}{E}}_{\beta_j})=1$,
$\chi({\stackrel{\circ}{E}}_{\alpha_j})=-1$,
and $\chi({\stackrel{\circ}{E}}_\sigma)=0$ for all other components
$E_\sigma$ of the exceptional divisor ${\cal D}$. Therefore
\begin{equation}\label{eq3}
\zeta_C(t)={\frac
{\prod\limits_{j=1}^s(1-t^{\beta_j})}
{\prod\limits_{j=0}^s(1-t^{\alpha_j})}
}.
\end{equation}
Comparing (\ref{eq1}) and (\ref{eq3}) and one has

\begin{theorem}\label{theo1}
The Poincar\'e series $P_C(t)$ of the ring of functions of an irreducible
plane curve singularity $C$ (or rather the rational function represented by it)
coincides with the zeta-function $\zeta_C(t)$ of the singularity $C$:
$P_C(t)=\zeta_C(t)$.
\end{theorem}

Up to now the authors do not know an explanation of this coincidence.
We consider it as an
exciting challenge.

\section{The extended semigroup of a plane curve singularity
and the Euler characteristic of its projectivization}\label{sec2}

It appears that one of possible generalizations of the result of
Theorem~\ref{theo1}
uses the notion of the extended semigroup of a curve singularity
(\cite{CDG1}). According to this generalization one can say that,
for an irreducible plane curve singularity $C$,
the coefficients of the Poincar\'e series $P_C(t)$ which are
equal to $0$ or $1$ are Euler characteristics of the empty set and of a
point respectively.

Let $(C, 0)\subset(\C^n, 0)$ be an arbitrary (reduced) curve
singularity (not necessarily plane) and let $C=\bigcup\limits_{i=1}^{r}C_i$
be its representation as the union of irreducible components $C_i$.
Let $\C_i$ be the complex line with the coordinate $\tau_i$ ($i=1, \ldots, r$)
and let $\varphi_i:(\C_i, 0)\to(\C^n, 0)$ be parameterizations
(uniformizations) of the branches $C_i$ of the curve $C$, i.e., germs of
analytic maps such that ${\rm{Im}}\,\varphi_i=C_i$ and $\varphi_i$ is an
isomorphism between $\C_i$ and $C_i$ outside of the origin. Let
${\cal O}_{\C^n, 0}$ be the ring of germs of holomorphic functions at the
origin in $\C^n$. For a germ $g\in{\cal O}_{\C^n, 0}$, let
$v_i=v_i(g)\in\Z_{\ge 0}$ and $a_i=a_i(g)\in\C^*$ be the degree of the
leading term and the coefficient at it in the power series
decomposition of the composition $g\circ\varphi_i:(\C_i,0)\to \C$:
$g\circ\varphi_i(\tau_i)=a_i\cdot\tau_i^{v_i}+{~terms~of~higher~degree}$,
$a_i\ne 0$ (if $g\circ\varphi_i(t)\equiv 0$, $v_i(g)=+\infty$
and $a_i(g)$ is not defined). The numbers $v_i(g)$ and $a_i(g)$ are defined
for elements $g$ of the ring ${\cal O}_C$ of functions on the curve $C$
as well. Let $v(g)=\Vert \vv(g) \Vert = v_1(g)+\ldots+v_r(g)$.

The semigroup $S_C$ of the curve singularity $C$ is the subsemigroup of
$\Z_{{\ge 0}}^r$ which consists of elements of the form
$\vv(g)=(v_1(g), \ldots, v_r(g))$ for all germs $g\in {\cal O}_C$ with
$v_i(g)<\infty$; $i=1, \ldots, r$. The set $\widehat S_C$ of elements
of the form $(\vv(g);\aa(g))=(v_1(g), \ldots, v_r(g); a_1(g), \ldots,
a_r(g)) \in \Z_{{\ge0}}^r\times(\C^*)^r$ for all germs $g\in {\cal O}_C$
with $v_i(g)<\infty$, $i=1, \ldots, r$, is a subsemigroup of
$\Z_{{\ge0}}^r\times(\C^*)^r$ and is called the extended semigroup of
the curve singularity $C=\bigcup\limits_{i=1}^{r}C_i$. The extended
semigroup $\widehat S_C$ is well defined (does not depend on the choice of
the parameterizations $\varphi_i$) up to a natural equivalence relation.

For plane curve singularities, the extended semigroup is not a topological
invariant, but reflects some moduli of it (see \cite{CDG1}).
Namely, it determines the exceptional divisor of the minimal embedded resolution
of a plane curve singularity up to projective equivalence, i.e., up to
projective equivalences of its irreducible components (each of them is
isomorphic to the projective line $\C\P^1$) which preserve all the
intersection points with other components of the total transform
$\pi^{-1}(C)$ of the curve $C$.

There is a natural semigroup homomorphism (projection) $p:\Z_{{\ge0}}^r\times
(\C^*)^r \to \Z_{{\ge0}}^r$ which sends $(\vv, \aa)$ to $\vv$. It maps the
extended semigroup $\widehat S_C$ of the curve $C$ onto its (usual) semigroup
$S_C$. The preimages $F_{\vv}=p^{-1}(\vv) \subset\{\vv\}\times(\C^*)^r \subset
\{\vv\}\times\C^r$ of the map $p:\widehat S_C\to \Z_{{\ge0}}^r$ are called
fibres of the extended semigroup $\widehat S_C$. $F_{\vv}$ is not empty if and
only if $\vv\in S_C$.

The fibres $F_{\vv}$ can be described in the following way.
For $\vv=(v_1, \ldots, v_r)\in \Z^r$ (not only for $\vv\in\Z^r_{\ge0}$),
let $J(\vv)$ be the ideal in the ring ${\cal O}_{\C^n, 0}$ of germs
of functions in $n$ variables
which consists of germs $g\in{\cal O}_{\C^n, 0}$ such that $\vv(g)\ge\vv$
(i.e., $v_i(g)\ge v_i$ for $i=1, \ldots, r$). There is a
natural linear map $j_{\vv}:J(\vv)\to\C^r$, which sends $g\in J(\vv)$ to
$\aa=(a_1, \ldots, a_r)\in\C^r$, where $a_i$ is the coefficient in the power
series decomposition $g\circ\varphi_i(\tau_i)=a_i\tau_i^{v_i}+
{~terms~of~higher~degree}$ (coefficients $a_i$ can be equal to zero).
The kernel of the map $j_{\vv}$ coincides with the ideal $J(\vv+\1)$
($\1=(1, \ldots, 1)$\,). Therefore the image ${\rm Im}\,j_\vv \subset\C^r$
of it is isomorphic to the vector space $C(\vv)=J(\vv)/J(\vv+{\1})$.
One can easily see that $F_{\vv}={\rm Im}\,j_\vv\cap(\C^*)^r$ (under the
natural identification of $\{\vv\}\times(\C^*)^r$ and $(\C^*)^r$). Therefore,
for $\vv\in S_C$, the fibre $F_{\vv}$ is the complement to an arrangement of
linear hyperplanes in a linear space. One can show that dimensions of the
fibres $F_{\vv}$ (or of the spaces $C(\vv)$\,) and combinatorial types of the
corresponding arrangements of hyperplanes ${\rm Im}\,j_\vv\cap
(\C^r\setminus(\C^*)^r)\subset{\rm Im}\,j_\vv$ are topological invariants
of a plane curve singularity (i.e., do not change along the stratum
$\mu=const$).

We shall be mostly interested in Euler characteristics of spaces under
consideration. The usual definition of the Euler characteristic of a
topological space (say, of a $CW$--complex) $X$ is $\chi(X)=
\sum\limits_{q\ge 0}(-1)^q\dim H_q(X;\R)$. If $X=X_1\cup X_2$ where the
spaces ($CW$--complexes) $X$, $X_1$ and $X_2$ are compact, one has
$\chi(X)=\chi(X_1)+\chi(X_2)-\chi(X_1\cap X_2)$. Therefore for compact
spaces the Euler characteristic possesses the additivity property.
This permits to consider it as a generalized (nonpositive) measure on
the algebra of such spaces. However spaces we are interested in (e.g.,
fibres of the extended semigroup of a curve) are noncompact semialgebraic
sets (complex or real). The Euler characteristic defined above does not
possess the additivity property for such spaces. For example, let $X$ be the
circle $S^1$, let $X_1$ be a point of $X$, and let $X_2=X\setminus X_1$ be
a (real) line. Then one has $\chi(X)=0$, $\chi(X_1)=\chi(X_2)=1$,
$\chi(X_1\cap X_2)=\chi(\emptyset)=0$, and $0\ne 1+1-0$.

In order to have the desired additivity property one should define the
Euler characteristic $\chi(X)$ of a space $X$ (say, semialgebraic) as
$$\sum\limits_{q\ge 0}(-1)^q\dim H_q(X^*, *;\R),$$ where $X^*$ is the
one-point compactification of the space $X$ (if $X$ is compact, the
one-point compactification of it is the disjoint union of $X$ with a
point), $*$ is the added ("infinite") point. We shall use this definition.
One can show that defined this way the Euler characteristic does possess
the additivity property for semialgebraic sets (in the example above
$\chi(X)=0$, $\chi(X_1)=1$, $\chi(X_2)=-1$). Moreover, a semialgebraic
set $X$ can be represented as a disjoint union of {\bf a finite number}
of open cells so that the boundary of a cell of some dimension (its
closure minus it itself) lies in the union of cells of smaller dimensions.
(This does not mean a representation of the space $X$ as a
$CW$--complex since in general (for noncompact sets) one does not have
maps of closed balls into $X$ which determine the cells. For example the
real line $\R^1$ is simply one cell of dimension 1.) One can see that the
Euler characteristic of the (semialgebraic) set $X$ is equal to the
alternative sum of numbers of cells of different dimensions.

By a graded topological space $X$ (with a $Z[[t_1, \ldots, t_r]]$--grading)
we shall mean a collection of topological spaces $X_\vv$ corresponding to
elements $\vv=(v_1, \ldots, v_r)$ of $\Z^r_{\ge 0}$. We shall (formally)
write $X=\sum\limits_{\vv\in\Z^r_{\ge 0}} X_\vv \cdot \tt^\vv$ where
$\tt=(t_1, \ldots, t_r)$, $\tt^\vv=t_1^{v_1}\cdot\ldots\cdot t_r^{v_r}$.
The (disjoint) union and the product of topological spaces defines the
sum and the product of graded spaces. By definition the Euler
characteristic of the graded space $X=\sum X_\vv \cdot \tt^\vv$ is the
power series $\chi(X)=\sum \chi(X_\vv) \cdot \tt^\vv$. One has
$\chi(X+Y)=\chi(X)+\chi(Y)$, $\chi(X\times Y)=\chi(X)\cdot\chi(Y)$. A
graded map from $X=\sum X_\vv \cdot \tt^\vv$ to $Y=\sum Y_\vv \cdot \tt^\vv$
is a collection of maps $X_\vv \to Y_\vv$ ($\vv\in\Z^r_{\ge 0}$).

The extended semigroup $\widehat S_C$ of a curve singularity $C=\bigcup
\limits_{i=1}^{r}C_i$ is in a natural sense a graded space:
$\widehat S_C=\sum F_\vv \cdot \tt^\vv$. The extended semigroup
$\widehat S_C$ (i.e., each fibre $F_\vv$) has a natural (free)
$\C^*$--action: multiplication of an element $\aa=(a_1, \ldots, a_r)\in
(\C^*)^r$ by a constant. Therefore $\chi(\widehat S_C)=0$. In order to
have a nontrivial Euler characteristic one can (or even should) "kill"
this $\C^*$. One way to do that is to factorize the extended semigroup
by this $\C^*$--action, i.e., to consider the projectivization of it.
(Another (somewhat dual) way was used in \cite{DL}).

Thus, let $\P\widehat S_C=\sum\P F_\vv\cdot t^\vv$ be the projectivization of the extended semigroup $\widehat S_C$ ($\P\widehat S_C \subset \Z_{\ge 0}
\times \P((\C^*)^r)$. This projectivization itself is a (graded) semigroup.
Its fibres $\P F_\vv$ are complements to arrangements of projective
hyperplanes in complex projective spaces. Let $X_C(t)$ be the Euler
characteristics of the projectivization $\P\widehat S_C$ of the extended
semigroup $\widehat S_C$, i.e., $X_C(t)= \sum \chi(\P F_\vv)\cdot t^\vv$.

It is not difficult to understand that, for a reducible plane curve
singularity ($r > 1$), the Euler characteristics $\chi(\P F_\vv)$ of
the fibres of the projectivization of the extended semigroup are equal
to zero for $\Vert\vv\Vert$ big enough and therefore $X_C(t)$ is in
fact a polynomial.

For an irreducible curve singularity $C$ ($r=1$), each fibre $F_\vv$
($\vv\in\Z^r_{\ge 0}$) is either empty (if $\vv\not\in S_C$) or isomorphic
to $\C^*$ (if $\vv\in S_C$). Thus its projectivization $\P F_\vv$ is empty
for $\vv\not\in S_C$ and is a point for $\vv\in S_C$ and
$$X_C(t)=\sum\limits_{\vv\in S_C}t^v=P_C(t).$$
Therefore, for a reducible curve singularity $C=\cup_{i=0}^r C_i$,
the polynomial $X_C(t)$ to some extend can be considered as a possible
generalization of the Poincar\'e series $P_C(t)$.

\section{Poincar\'e series of the ring of functions of a reducible
curve singularity}\label{sec3}

Another (by definition, not by essence: see Theorem~\ref{theo2})
generalization of the Poincar\'e series of the ring of of germs of functions
for a reducible curve singularity $C=\cup_{i=0}^r C_i\subset(\C^n, 0)$
can be constructed in the following way. For $\vv \in \Z^r$ (not only
for $\vv \in \Z_{\ge 0}^r$), let $c(\vv)$ be the (non-negative) integers
$\dim J(\vv)/J(\vv+\1)$. Let ${\cal L}(t_1,\ldots,t_r) =
\sum\limits_{\vv \in {\Z}^r} c(\vv) t^\vv$.

For $r>1$, the integers $c(\vv)$ can be positive for $\vv$ with (some)
negative
components $v_i$ as well. For example, if there exists a germ
$g\in {\cal O}_{\C^n, 0}$ with $v_1(g)=v_1^*$, then for any $v_2$,
\dots, $v_r$ such that $v_i\le v_i(g)$ (including negative ones), the
germ $g$ represents a non-trivial element in the factor $J(\vv)/J(\vv+\1)$
where $\vv=(v_1^*, v_2, \ldots, v_r)$, $J(\vv)=\{g\in{\cal O}_{\C^n, 0}:
\vv(g)\ge\vv\}$, and therefore $c(\vv)>0$. Therefore ${\cal L}(t_1,\ldots,t_r)$ is not a
power series but an element of the set $\Z[[t_1, \ldots, t_r, t_1^{-1},
\ldots, t_r^{-1}]]$ of Laurent series in $t_1, \ldots, t_r$ infinite in all
the directions. The set $\Z[[t_1, \ldots, t_r,t_1^{-1}, \ldots, t_r^{-1}]]$
is not a ring (generally speaking, the product of its elements is not
defined), but it is a $\Z[t_1, \ldots, t_r]$-- or a
$\Z[t_1, \ldots, t_r,t_1^{-1}, \ldots, t_r^{-1}]$--module. One can understand
that along each line in the lattice $\Z^r$ parallel to a coordinate one
the coefficients $c(\vv)$ stabilize in each direction. This is the reason why
$P'_C(t_1, \cdots ,t_r)= {\cal L}(t_1, \cdots, t_r)\cdot
\prod\limits_{i=1}^r (t_i-1)$ is a polynomial. Moreover, one can show that,
for $r>1$, the polynomial $P'_C(t_1, \ldots, t_r)$
is divisible by $(t_1\cdot\ldots\cdot t_r-1)$. Let $P_C(t_1,\ldots,t_r)=
P'_C(t_1,\ldots,t_r)/(t_1\cdot\ldots\cdot t_r-1)$. Then $P_C(\tt)$ is a
polynomial for $r>1$ and is a power series for $r=1$ (i.e. for
an irreducible curve singularity). In the last case $P_C(\tt)$ coincides
with the Poincar\'e series of the ring of functions on the curve discussed
above and, thus, it can be considered as a generalization of it. We shall
call $P_C(t_1, \ldots, t_r)$ the Poincar\'e series of the curve singularity
$C=\cup_{i=1}^r C_i$.

Though the definitions of the polynomials (series for $r=1$) $X_C(t)$
and $P_C(t)$ look rather different, in fact they coincide (for any
curve, not only for a plane one).

\begin{theorem}\label{theo2} 
For a curve singularity $C=\bigcup\limits_{i=1}^r C_i \subset (\C^n, 0)$
one has $$X_C(t)=P_C(t).$$
\end{theorem}

The {\bf proof} consists in application of the inclusion--exclusion
formula for computing the Euler characteristic of the projectivization
$\P F_\vv$ as the complement to an arrangement of projective hyperplanes
in a projective space (taking into account that $\chi(\P(\C^n))=
\chi(\CP^{n-1})=n$).

\section{The Alexander polynomial of a plane curve singularity and the
Poincar\'e polynomial}
\label{sec4}

The Alexander polynomial (in $r$ variables) is an invariant of a link
with $r$ (numbered) components in the sphere $S^3$. The general definition
can be found, e.g., in \cite{EN}. To a plane curve singularity
$C=\bigcup\limits_{i=1}^r C_i\subset (\C^2, 0)$ there corresponds
the link $C\cap S^3_\varepsilon$ in the 3--sphere $S^3_\varepsilon$
of radius $\varepsilon$ centred at the origin in the complex plane $\C^2$
with $\varepsilon$ small enough. Let $\Delta^C(t_1, \ldots, t_r)$ be
the Alexander polynomial of this link ($\equiv$ the Alexander polynomial
of the curve $C$). For such a link (an algebraic one) we rather use
not the general definition of the Alexander polynomial, but a formula for
it in terms of an embedded resolution $\pi:({\cal X}, {\cal D})\to(\C^2, 0)$
of the curve singularity $C$.

Let the curve $C=\bigcup\limits_{i=1}^r C_i$ be given by an equation
$f=0$ ($f:(\C^2, 0)\to(\C, 0)$) and let $f=\prod\limits_{i=1}^r f_i$
where $f_i=0$ are equations of the components $C_i$
of the curve $C$ ($i=1, \ldots, r$).
For each component $E_\sigma$ of the exceptional divisor ${\cal D}=
\pi^{-1}(0)$ let $m_i^\sigma$ be the multiplicity of the lifting
$f_i\circ\pi$ of the function $f_i$ to the space ${\cal X}$ of the
resolution along the component $E_\sigma$, ${\mm}^\sigma=
(m_1^\sigma, \ldots, m_r^\sigma)$. One has $\sum\limits_{i=1}^r m_i^\sigma=
m^\sigma$. As above let ${\stackrel{\circ}{E}}_\sigma$ be the smooth part
of the component $E_\sigma$, i.e., $E_\sigma$ minus its intersections
with all other components of the total transform $\pi^{-1}(C)$ of the
curve $C$. These intersection points on the component $E_\sigma$ are
in one-to-one correspondence with connected components of the complement
$\mbox{$(f\circ\pi)^{-1}(0)\setminus {\stackrel{\circ}{E}}_\sigma$}$. We
shall say that an intersection point is {\bf essential} if the corresponding
connected component contains the strict transform of a component of the
curve $C$. Let $s_\sigma$ be the number of essential points on the component
$E_\sigma$, and let ${\widetilde E}_\sigma$ be the complement to the set
of essential points in $E_\sigma$.

\begin{theorem}\label{theo3} {\rm (\cite{EN})} For $r>1$,
\begin{equation}\label{eq4}
\Delta^C(t_1, \ldots, t_r)=\prod\limits_{E_\sigma\subset{\cal D}}
\left(1-\tt^{{\mm}^\sigma}\right)^{-\chi({\stackrel{\circ}{E}}_\sigma)}.
\end{equation}
\end{theorem}

One can see that this formula is an analogue of the formula (\ref{eq2})
for the zeta--function $\zeta_C(t)$ of the monodromy transformation
of the curve $C$. Moreover one has $\zeta_C(t)=\Delta^C(t, \ldots, t)$
(for $r>1$).

\begin{remarks}
1. For an irreducible curve singularity (i.e., for $r=1$) the right-hand side of the
formula (\ref{eq4}) is not a polynomial, but a polynomial (in fact the
Alexander polynomial in one variable) divided by $(1-t)$.\newline
2. Generally speaking the Alexander polynomial is defined only up
to multiplication by a monomial $\pm{\underline{t}}^{\underline{k}}$,
$\kk=(k_1, \ldots, k_r)\in\Z^r$.
For algebraic links the formula (\ref{eq4}) fixes the choice of the
Alexander polynomial in such a way that it is really a polynomial
(i.e., does not contain monomials with negative powers) and its value
at the origin ($\tt=0$) is equal to 1.
\end{remarks}

\begin{theorem}\label{theo4}
For a plane curve singularity $C=\bigcup_{i=1}^r C_i\subset (\C^2, 0)$,
$r>1$,
$$
\Delta^C(\tt)=X_C(\tt)=P_C(\tt).
$$
\end{theorem}

The main course of the {\bf proof} goes as follows. First we construct a
graded space (in fact a graded semigroup) $Y$ such that its Euler
characteristic is equal to the Alexander polynomial $\Delta^C(t_1, \ldots,
t_r)$ together with a map (a semigroup homomorphism) to the projectivization
$\P\widehat{S}_C$ of the extended semigroup $\widehat{S}_C$. Let
$\pi:({\cal X}, {\cal D})\to(\C^2, 0)$ be an embedded resolution of the
curve singularity $C\subset(\C^2, 0)$, let the exceptional divisor ${\cal D}$
of the resolution $\pi$ be the union of irreducible components $E_\sigma$
($\sigma\in\Gamma$). For a topological space $X$, let $S^kX=X^k/S_k$
($k\ge 0$) be the $k$th symmetric power of the space $X$, i.e., the space
of subsets of the space $X$ with $k$ elements ($S^0X$ is a point).

Let
$$
Y=\sum\limits_{\{k_\sigma\}}\left(\prod\limits_\sigma S^{k_\sigma}
{\stackrel{\circ}{E}}_\sigma\cdot\tt^{k_\sigma{\mm}^\sigma}\right),
$$
where the union (sum) is over all sets of nonnegative integers $k_\sigma$,
$\sigma\in\Gamma$, $S^k{{\stackrel{\circ}{E}}_\sigma}$ is the $k$-th
symmetric power of the nonsingular part ${\stackrel{\circ}{E}}_\sigma$
of the component $E_\sigma$ of the exceptional divisor. The graded
topological space $Y$ is a (graded) semigroup with respect to the
operation defined by the union of subsets.

There is a natural map (a graded semigroup homomorphism)
$\Pi:Y\to\P\widehat S_C$ defined as follows. An element
$$
y\in Y=\prod\limits_{\sigma}\left(\bullet+
S^1{\stackrel{\circ}{E}}_\sigma\cdot\tt^{{\mm}^\sigma}
+S^2{\stackrel{\circ}{E}}_\sigma\cdot\tt^{2{\mm}^\sigma}+\ldots\right),
$$
where $\bullet$ is a point ($=S^0{\stackrel{\circ}{E}}_\sigma$),
is represented by a subset of "nonsingular" points of the exceptional divisor
$\cal D$ (i.e., of ${\stackrel{\circ}{\cal D}}=\bigcup\limits_\sigma
{\stackrel{\circ}{E}}_\sigma$) with $k_\sigma$ points $P_1^\sigma$,
\dots, $P_{k_\sigma}^\sigma$ on the component ${\stackrel{\circ}{E}}_\sigma$.
For a point $A\in {\stackrel{\circ}{\cal D}}$, let
$\widetilde L_A$ be the germ of a nonsingular (complex analytic)
curve transversal to the exceptional divisor ${\cal D}$ at the point
$A$. Let the image $L_A=\pi(\widetilde L_A)\subset(\C^2, 0)$ of the curve
$\widetilde L_A$ be given by an equation $g_A=0$,
$g_A\in{\cal O}_{\C^2, 0}$. Let $g_y=\prod\limits_\sigma
\prod\limits_{j=1}^{k_\sigma}g_{P_j^\sigma}$. Then
$$
\Pi(y):=(v_1(g_y),\ldots,v_r(g_y);a_1(g_y):\ldots:a_r(g_y)).
$$

To prove that the map $\Pi$ is well defined, i.e., that
$\Pi(y)\in\P\widehat S_C$ does not depend on the choice of curves
${\widetilde L}_A$ and of the equations $g_A=0$, suppose that
${\widetilde L}_A^\prime$ is another germ of a nonsingular (complex
analytic) curve transversal to the exceptional divisor $\cal D$
at the point $A\in{\stackrel{\circ}{D}}$, $g_A^\prime=0$ is an equation
of the curve $L_A^\prime=\pi({\widetilde L}_A^\prime)\subset(\C^2, 0)$. Let
$g_y^\prime=
\prod\limits_\sigma\prod\limits_{j=1}^{k_\sigma}g_{P_j^\sigma}^\prime$,
let $\widetilde g_y=g_y\circ\pi$ and
$\widetilde g_y^\prime=g_y^\prime\circ\pi$ be the liftings of the germs
$g_y$ and $g_y^\prime$ to the space $\cal X$ of the resolution, and let
$\psi=\widetilde g_y^\prime/\widetilde g_y$ be their ratio. The meromorphic
function $\psi$ on $\cal X$ has (simple) zeros along the curves
${\widetilde L}_{P_j^\sigma}^\prime$ and (also simple) poles along the
curves ${\widetilde L}_{P_j^\sigma}$. Therefore the restriction of the
function $\psi$ to the exceptional divisor $\cal D$ has neither zeroes
no poles, i.e., it is a regular (holomorphic) function without zeroes
and thus it is a constant (say, $c$) on $\cal D$. It implies that
$\vv(g_y^\prime)=\vv(g_y)$, $\aa(g_y^\prime)=c\cdot\aa(g_y)$ and therefore
the points $(\vv(g_y^\prime);a_1(g_y^\prime):\ldots:a_r(g_y^\prime))$ and
$(\vv(g_y);a_1(g_y):\ldots:a_r(g_y))$ in the projectivization
$\P\widehat S_C$ coincide.

To compute the Euler characteristic $\chi(Y)$ of the graded space $Y$
one uses the following statement.

\begin{lemma}
For a topological space $X$,
$$
\chi(\bullet+S^1X\cdot t+S^2X\cdot t^2+\cdots)=(1-t)^{-\chi(C)}.
$$
\end{lemma}

\begin{corollary}
$\displaystyle{
\chi(Y)=\prod\limits_\sigma
\left(1-\tt^{{\mm}^\sigma}\right)^{-\chi({\stackrel{\circ}{E}}_\sigma)}=
\Delta^C(t_1,\ldots, t_r).
}
$
\end{corollary}

Let $\underline{V}$ be an arbitrary point of the lattice $\Z_{\ge0}^r$
and suppose that the resolution $\pi:({\cal X}, {\cal D})\to(\C^2, 0)$
of the curve $C$ is such that, for any germ $g\in{\cal O}_{\C^2,0}$
with $\vv(g)\le \underline{V}$, the strict transform
$\overline{(g\circ\pi)^{-1}(0)\setminus{\cal D}}$ of the curve $\{g=0\}$
intersects the exceptional divisor $\cal D$ only at nonsingular points of it.
One can get such a resolution from the minimal one making sufficiently many
blow-ups at intersection points of components of the total transform
$\pi^{-1}(C)$ of the curve $C$. It is not difficult to see that in this case,
for $\vv\le\underline{V}$, the image $\Pi(Y_\vv\cdot\tt^\vv)$ of the part
of grading $\vv$ of the space $Y$ coincides with the projectivization
$\P F_\vv$ of the fibre $F_\vv$ of the extended semigroup $\widehat S_C$
(or rather with the corresponding summand $\P F_\vv$). Therefore in order
to prove that $\chi(\P\widehat S_C)= \Delta^C(t_1, \ldots, t_r)$ it is
sufficient to show that $\chi({\mbox{Im\,}}\Pi)=\chi(Y)$.
There would be no problem to show this if the map $\Pi$ would be injective.
This is not the case. Thus one has to analyse 
how non-injective is it.

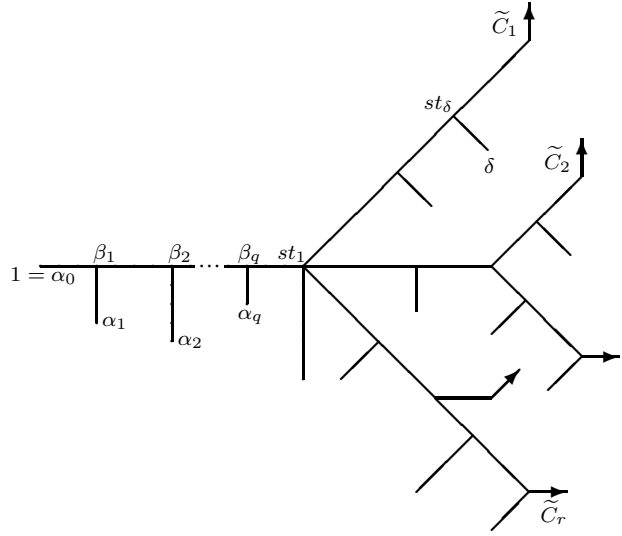
\begin{figure}
$$
\unitlength=0.50mm
\begin{picture}(120.00,100.00)(10,-40)
\thicklines
\put(-10,30){\line(1,0){41}}
\put(39,30){\line(1,0){21}}
\put(60,30){\line(1,0){50}}
\put(90,30){\line(0,-1){12}}

\put(60,30){\line(1,1){60}}
\put(120,90){\vector(0,1){10}}
\put(110,92){{\scriptsize $\widetilde C_1$}}
\put(85,55){\line(1,-1){9}}
\put(100,70){\line(1,-1){9}}
\put(108,55){{\scriptsize $\delta$}}
\put(92,72){{\scriptsize $st_\delta$}}

\put(60,30){\line(1,-1){60}}
\put(120,-30){\vector(1,0){10}}
\put(123,-38){{\scriptsize $\widetilde C_r$}}
\put(80,10){\line(-1,-1){10}}
\put(105,-15){\line(-1,-1){15}}
\put(120,-30){\line(-1,-1){10}}
\put(95,-5){\line(1,0){15}}
\put(110,-5){\vector(1,1){7.5}}

\put(110,30){\line(1,1){24}}
\put(122,42){\line(1,-1){9}}
\put(134,54){\vector(0,1){10}}
\put(124,56){{\scriptsize $\widetilde C_2$}}

\put(110,30){\line(1,-1){24}}
\put(119,21){\line(-1,-1){9}}
\put(134,6){\line(-1,-1){9}}
\put(134,6){\vector(1,0){10}}

\put(33,30){\circle*{0.5}}
\put(35,30){\circle*{0.5}}
\put(37,30){\circle*{0.5}}
\put(25,10){\line(0,1){20}}
\put(45,20){\line(0,1){10}}
\put(60,0){\line(0,1){30}}
\put(5,15){\line(0,1){15}}
\thinlines
\put(15,30){\circle*{1}}
\put(25,30){\circle*{1}}
\put(45,30){\circle*{1}}
\put(60,30){\circle*{1}}
\put(25,20){\circle*{1}}
\put(60,20){\circle*{1}}
\put(60,10){\circle*{1}}
\put(5,30){\circle*{1}}
\put(25,10){\circle*{1}}
\put(45,20){\circle*{1}}
\put(60,0){\circle*{1}}

\put(-10,30){\circle*{1}}
\put(-5,30){\circle*{1}}
\put(0,30){\circle*{1}}
\put(10,30){\circle*{1}}
\put(20,30){\circle*{1}}
\put(30,30){\circle*{1}}
\put(40,30){\circle*{1}}
\put(50,30){\circle*{1}}
\put(55,30){\circle*{1}}
\put(60,25){\circle*{1}}
\put(60,15){\circle*{1}}
\put(60,5){\circle*{1}}
\put(5,25){\circle*{1}}
\put(5,20){\circle*{1}}
\put(5,15){\circle*{1}}
\put(25,25){\circle*{1}}
\put(25,15){\circle*{1}}
\put(-18,26){{\scriptsize $1=\alpha_0$}}
\put(6.5,14){{\scriptsize$\alpha_1$}}
\put(26.5,9){{\scriptsize$\alpha_2$}}
\put(42.5,16){{\scriptsize$\alpha_q$}}
\put(4,32){{\scriptsize$\beta_1$}}
\put(24,32){{\scriptsize$\beta_2$}}
\put(42.5,32){\scriptsize{$\beta_q$}}
\put(53,32){\scriptsize{$st_1$}}
\end{picture}
$$
\caption{The dual resolution graph $\Gamma$ of a redicible curve.}
\label{fig2}
\end{figure}

One possibility to decrease the set where the map $\Pi$ is non-injective
is is to reduce the graded space $Y$ in the following way. Let $D^\prime$
be the union of components $E_\sigma$ of the exceptional divisor with
at least two essential points, i.e., with $s_\sigma\ge 2$, and let
$\Delta^\prime$ be the set of the corresponding vertices of the dual graph
of the resolution. Connected components of the complement
$D\setminus D^\prime$, which do not contain the starting divisor
${\bf 1}$, are tails of the dual graph $\Gamma$ of the
resolution and correspond to (some) dead ends $\delta$ of the
graph $\Gamma$. Let $\Delta$ be the set of these dead ends. For
$\delta\in\Delta$, let $st_\delta$ be the vertex of $\Delta^\prime$
such that $E_{st_\delta}$ intersects the corresponding connected
component of $\overline{D\setminus D^\prime}$ (see Fig.{\ref{fig2}}). Let
$$
Y^\prime=\prod\limits_{\sigma\in\Delta^\prime}(\bullet+
S^1{{\widetilde{E}}_\sigma}\cdot\tt^{{\mm}^\sigma}
+S^2{{\widetilde{E}}_\sigma}\cdot\tt^{2{\mm}^\sigma}+\ldots).
$$
Pay attention that the spaces ${\stackrel{\circ}{E}}_\sigma$
("nonsingular parts" of the components of the exceptional divisor)
in the definition of the graded space $Y$ are substituted here by the spaces
${\widetilde E}_\sigma$ (the components of the exceptional divisor
minus essential points of them).

For a dead end $\delta\in\Delta$, ${\mm}^{st_\delta}$ is a multiple of
${\mm}^{\delta}$: ${\mm}^{st_\delta}=(n_\delta+1)\cdot {\mm}^{\delta}$. Let
$Y_\delta=\sum\limits_{k=0}^{n_\delta}\bullet\cdot\tt^{k{\mm}^\delta}$.

Let $st_1$ be the first vertex with at least two essential points (see Fig.{\ref{fig2}}).
Let us assume that ${\bf 1} \neq st_1$. Let $\beta_0={\bf 1}$, $\beta_1$,
\dots, $\beta_q$ be the dead ends of the graph $\Gamma$ which do not
belong to $\Delta$, and let $\alpha_j$ ($j=1, \ldots, q$) be the star
point of the graph $\Gamma$ which corresponds to the dead end $\beta_j$,
$\alpha_1 < \alpha_2 < \ldots < \alpha_q$. Let $S_1$ be the subsemigroup
of the semigroup $S_C$ generated by the multiplicities ${\mm}^{\beta_0}$,
${\mm}^{\beta_1}$, \dots, ${\mm}^{\beta_q}$.
$S_1$ coincides with the semigroup generated by all
the multiplicities ${\mm}^\sigma$ with $E_\sigma$ from the connected
component of $\overline{D\setminus D^\prime}$ which contains the starting
divisor ${\bf 1}$ and is similar to a subsemigroup
of $\Z_{\ge 0}$ because it is contained in the line $L$ in
$\R^r\supset\Z_{\ge 0}^r$ which goes through the origin and
the point ${\mm}^{\beta_0}$. Let $S_1^\prime$ be the subset of
the semigroup $S_1$ which consists of elements
$\mm\in S_1$ such that ${\mm}-{\mm}^{st_1}\not\in S_1$.
Let $Y_1=\sum\limits_{{\mm}\in S_1^\prime}\bullet\cdot\tt^{\mm}$
($S^\prime_1$ is a finite set).
If ${\bf 1} = st_1$ we simply put $Y_1 = \bullet$.

Let $\displaystyle{\widetilde Y =
Y^\prime\times Y_1\times\prod\limits_{\delta\in\Delta} Y_\delta}$.
One can show that $\chi{(\widetilde{Y})} = \chi(Y)$.

There exists a map $\widetilde{\Pi}:
\widetilde{Y}\to \P \widehat{S}_{C}$ such that $Im\,\widetilde{\Pi}
=Im\,\Pi$.
To define it, one can say that $\widetilde Y$ is a subset of the
graded semigroup
$$\widetilde Y^*=Y^\prime\times\left(\sum\limits_{{\mm}\in S_1}
\bullet\cdot\tt^{\mm}\right)\times
\prod\limits_{\delta\in\Delta}\left(
\sum\limits_{k=0}^{\infty}\bullet\cdot\tt^{k{\mm}^\delta}\right)
$$
(each factor of $\widetilde Y^*$ is a graded semigroup)
and the map $\widetilde\Pi$ is the restriction of a graded semigroup
homomorphism $\widetilde Y^*\to\P\widehat S_C$. Because of that it
should be defined for points of
$\bigcup\limits_{\sigma\in\Delta^\prime}\widetilde E_\sigma$ and also for
"monomials" of the form $\bullet\cdot\tt^{{\mm}^\delta}$ for
$\delta\in\Delta$ and $\bullet\cdot\tt^{{\mm}^{\beta_i}}$ for
$i=0, 1, \ldots, q$.

For a point $A$ of ${\stackrel{\circ}{E}}_\sigma$,
$\sigma\in\Delta^\prime$, (or rather for the monomial
$[A]\cdot\tt^{{\mm}^\sigma}$) $\widetilde\Pi$ coincides with $\Pi$.
A point of $\widetilde E_\sigma\setminus{\stackrel{\circ}{E}}_\sigma$,
$\sigma\in\Delta^\prime$, corresponds either to a dead
end $\delta\in\Delta$ (and in this case $\sigma=st_\delta$) or to the
initial divisor ${\bf 1}$ (in this case $\sigma=st_{1}$). In the
first case one puts $\widetilde\Pi([A]\cdot\tt^{{\mm}^\sigma})=
(n_\delta+1)\cdot\Pi([A_\delta]\cdot\tt^{{\mm}^\delta})$ for any
point $A_\delta\in{\stackrel{\circ}{E}}_\delta$ (see the definition
of $n_\delta$ above); in the second case there exists
$\ell_0$, \ldots, $\ell_q$ such that $\mm^\sigma=\mm^{st_1} =
\sum\limits_{i=0}^q  \ell_i\mm^{\beta_i}$ and one puts
$\widetilde\Pi([A]\cdot\tt^{{\mm}^\sigma})=
\sum\limits_{i=0}^q\ell_i\cdot\Pi([A_{\beta_i}]\cdot\tt^{{\mm}^{\beta_i}})$
for any points $A_{\beta_i}\in{\stackrel{\circ}{E}}_{\beta_i}$.
One puts $\widetilde\Pi(\bullet\cdot\tt^{{\mm}^\delta})=
\Pi([A_\delta]\cdot\tt^{{\mm}^\delta})$ for any point
$A_\delta\in{\stackrel{\circ}{E}}_\delta$, $\delta\in\Delta$,
$\widetilde\Pi(\bullet\cdot\tt^{{\mm}^{\beta_i}})=
\Pi([A_{\beta_i}]\cdot\tt^{{\mm}^{\beta_i}})$ for any point
$A_{\beta_i}\in{\stackrel{\circ}{E}}_{\beta_i}$, $i=0, 1, \ldots, q$
(one can easily see that the result does not depend on the choice of the
points $A_\delta$, $A_{\beta_i}$ in these cases).

It is not difficult to see that $Im\,\widetilde{\Pi}=Im\,\Pi$. Now
Theorem~\ref{theo4} follows from the statement that $\chi(\widetilde{Y})=
\chi(\mbox{Im }\widetilde{\Pi})$. To prove this, one analyses places where the map $\widetilde{\Pi}$ is not injective. At such places we indicate some parts of the space $\widetilde Y$
which are fibred into complex tori $(\C^*)^{s-1}$, $s\ge 2$, (and thus
have zero Euler characteristic) and which can be removed without
changing the image $\widetilde{\Pi}(\widetilde Y)$.
To obtain such a proof, the following three types of statements are used.

\smallskip
{\bf 1.} Connected components of the space $\widetilde Y$
are labelled by multi-indices
$$
\kk = \left\{
\{k_{\sigma}\}_{\sigma\in \Delta'}, \mm,
\{k_{\delta}\}_{\delta\in \Delta} \right\}
$$
with $k_{\sigma}\ge 0$ for each  $\sigma\in \Delta'$, $\mm\in S_1^\prime$,
and $0\le k_{\delta}\le n_{\delta}$ for each $\delta\in \Delta$.
A vertex $\sigma\in \Delta'$ is called a cut of a multi-index
$\kk=\{k_\sigma, \mm, k_\delta\}$
if $k_\sigma\ge s_\sigma$. One shows that if $\widetilde\Pi(y_1)=
\widetilde\Pi(y_2)$ for two different elements $y_1$ and $y_2$ from
$\widetilde Y$ corresponding to multi-indices $\kk^i$, $i=1, 2$
(i.e., at a place where the map $\widetilde\Pi$ is not injective),
then each multi-index $\kk^i$ ($i=1, 2$) has a cut. Moreover
one proves more fine statements about the distribution of the cuts
of the multi-indices $\kk^1$ and $\kk^2$ on the tree $\Delta^\prime$.
Namely, up to the numbering of $y_1$ and $y_2$ there exists a cut
$\sigma$ of the multi-index $\kk^1$ such that for each strict transform
$\widetilde C_i$ of a branch $C_i$ of the curve $C$ greater than $\sigma$
there is a cut of the multi-index $\kk^2$ on the geodesic from $\sigma$
to $\widetilde C_i$. (For vertices $\sigma_1$ and $\sigma_2$ of the
dual graph $\Gamma$ of the resolution, one says that $\sigma_1\ge\sigma_2$
if $\sigma_2$ lies on the geodesic between $\sigma_1$ and the starting
divisor $\bf 1$.) This is the most complicated (combinatorial) part of
the proof.

\smallskip
{\bf 2.} It is not difficult to describe the difference between the images
$\widetilde\Pi(y)$ and $\widetilde\Pi(\widehat{y})$ for points $y$ and
$\widehat y$ from the space $\widetilde Y$ corresponding to one and the
same multi-index $\kk$ and such that there exists only one vertex
$\sigma\in\Delta^\prime$ for which the sets of points on
${\stackrel{\circ}{E}}_{\sigma}$ for the elements $y$ and $\widehat y$
($P^\sigma_1$, \dots, $P^\sigma_{k_\sigma}$ and $\widehat{P}^\sigma_1$,
\dots, $\widehat{P}^\sigma_{k_\sigma}$) respectively) are different.
Let $Q^\sigma_0$, $Q^\sigma_1$, \dots, $Q^\sigma_{s-1}$ be the essential
points on the component $E_\sigma$ of the exceptional divisor. Let
$\Psi$ be a meromorphic function on $E_\sigma$ ($\cong \C\P^1$) with
zeroes at the points $\widehat{P}^\sigma_1$, \dots,
$\widehat{P}^\sigma_{k_\sigma}$ and poles at the points $P^\sigma_1$,
\dots, $P^\sigma_{k_\sigma}$. Such function is well-defined up to
multiplication by a non-zero constant. For $i=1, \ldots, r$, let
$j(i)$ be defined by the condition that the connected component of
$\overline{\pi^{-1}(C)\setminus E_\sigma}$ which contains the strict
transform $\widetilde C_i$ of the branch $C_i$ of the curve $C$ intersects
the component $E_\sigma$ of the exceptional divisor at the point
$Q^\sigma_{j(i)}$. Then $\widetilde\Pi(\widehat y)$ is obtained from
$\widetilde\Pi(y)$ by multiplying the coordinate $a_i$ ($i=1, \ldots, r$)
by the constant $\Psi(Q^\sigma_{j(i)})$ (the value of the described
meromorphic function $\Psi$ at the corresponding essential point).
Therefore it is necessary to describe sets of possible values of the
function $\Psi$ at the essential points $Q^\sigma_0$, $Q^\sigma_1$,
\dots, $Q^\sigma_{s-1}$ for different sets of points $\widehat{P}^\sigma_1$,
\dots, $\widehat{P}^\sigma_{k_\sigma}$.

\smallskip
{\bf 3.} Let $Q_0$, $Q_1$, \dots, $Q_{s-1}$ be (different) points of a
projective line $E$, $\widetilde{E} = E- \{Q_\ell : \ell=0, 1,
\ldots, s-1\}$, let $P^o_1$, \dots, $P^o_k$ be $k$ points (not
necessarily different), different from $Q_0$, $Q_1$, \dots, $Q_{s-1}$.
Let $\Phi$ be the map from the symmetric power $S^k \widetilde{E}$ of
the space $\widetilde{E}$ to $\P ((\C^*)^s)$ defined in the following way.
For an element from $S^k\widetilde{E}$; i.e., for $k$ points $P_1$, \dots,
$P_k$ from $\widetilde{E}$, let $\psi$ be a meromorphic function on $E$
with zeroes at the points $P_1$, \dots, $P_k$ and poles at the points
$P^o_1$, \dots, $P^o_k$; let $\Phi(\{P_j\}) := (\psi(Q_0): \psi(Q_1):
\ldots : \psi(Q_{s-1}))$. Then, if $k\ge s-1$, one has $\mbox{Im }\Phi =
\P((\C^*)^s)$; if $k\le s-1$, $\Phi$ is an embedding. Moreover in both
cases the map $\Phi$ is a (smooth) locally trivial (in fact a trivial)
fibration over its image the fibre of which is a (complex) affine space
of dimension $\max(0, k-s+1)$.

In order to prove this one can proceed as follows. Without loss of
generality one can suppose that $P^o_1=P^o_2=\cdots = P^o_k=P^o$. Let us
choose an affine coordinate on the projective line $E$ such that
$P^o=\infty$. Then $\psi$ is a polynomial of degree $\le k$ with zeroes at
those points $P_1$, \dots, $P_k$ which are different from $P^o$. Let
$z_\ell$ be the coordinate of the point $Q_\ell$, $\ell=0, 1, \ldots, s-1$.

For $k\ge s-1$, the statement that the map $\Phi$ is onto can be reduced
to the following obvious one: for an arbitrary prescribed set of values
$\{\psi_0$, $\psi_1$, \dots, $\psi_{s-1}\}$, there exists a
polynomial $\psi$ of degree $\le k$ such that
$\psi(z_\ell)=\psi_\ell$, $\ell=0, 1, \ldots, s-1$. The statement that
$\Phi$ is a locally trivial fibration over its image
follows from the fact that if $\psi_1$ and $\psi_2$ are polynomials
with coinciding values at the points $Q_\ell$, $\ell=0, 1, \ldots, s-1$,
then $\psi_1=\psi_2+q(z)(z-z_0)(z-z_1)\cdot\ldots\cdot(z-z_{s-1})$
where $q(z)$ is an arbitrary polynomial of degree $k-s$.

For $k\le s-1$, the statement follows from the fact that such a
polynomial of degree $\le s-1$ is unique.

This statement implies that at places where the map $\widetilde\Pi$ is not
injective (and thus in the presence of cuts) one can find parts of the
graded space $\widetilde Y$ which have zero Euler characteristic and which
can be deleted without changing the image of the map $\widetilde\Pi$.

\smallskip
Combination of all these arguments (and/or of their versions) proves
Theorem~\ref{theo4}.

\section{The Alexander polynomial as an integral with respect to
the Euler characteristic}
\label{sec5}

It appears that Theorem~\ref{theo4} can be reformulated so that the
Alexander polynomial $\Delta^C(t_1, \cdots, t_r)$ will be written as
a certain integral with respect to the Euler characteristic over the
space ${\cal O}_{\C^2, 0}$ of germs of functions on the plane $\C^2$ at
the origin or rather over its projectivization $\P{\cal O}_{\C^2, 0}$.
In order to do that one has to define the notion of the Euler
characteristic for (some) subsets of the space $\P{\cal O}_{\C^2, 0}$.
The idea is the same which is used to define the notion of the Euler
characteristic (or of the generalized Euler characteristic with
values in the Grothendieck ring of semialgebraic sets localized by
the class of the complex line) for subsets of the space of arcs (in the
framework of the theory of the motivic integration; see, e.g., \cite{C}).
The definitions of the Euler characteristic of (some) subsets and
of the integral with respect to the Euler characteristic is practically
the same for the space ${\cal O}_{\C^n, 0}$ of germs of functions at
the origin in the space $\C^n$ and for its projectivization $\P{\cal O}_{\C^n, 0}$.
Since here we shall use the last one, we shall give them in this
setting.

Let $J^k_{\C^n,0}$ be the space of $k$-jets of functions at the origin in
$(\C^n,0)$ ($J^k_{\C^n,0}={\cal O}_{\C^n,0}/{\goth m}^{k+1}\cong
\C^{{n+k\choose k}}$, where $\goth m$ is the maximal ideal in the ring
${\cal O}_{\C^n,0}$). For a complex vector
space $L$ (finite or infinite dimensional) let $\P L=(L\setminus\{0\})/\C^*$
be its projectivization, let $\P^* L$ be the disjoint union of $\P L$
with a point (in some sense $\P^* L=L/\C^*$). One has natural maps
$\pi_k: \P{\cal O}_{\C^n,0} \to \P^* J^k_{\C^n,0}$ and
$\pi_{k,\ell}: \P^* J^k_{\C^n,0} \to \P^* J^\ell_{\C^n,0}$ for $k \ge \ell$.
Over $\P J^\ell_{\C^n,0} \subset \P^* J^\ell_{\C^n,0}$ the map $\pi_{k,\ell}$
is a locally trivial fibration, the fibre of which
is a complex vector space of dimension ${n+k\choose k}-{n+\ell\choose \ell}$.

\begin{definition}
A subset $X\subset \P{\cal O}_{\C^n,0}$ is said to be cylindric if
$X=\pi_k^{-1}(Y)$ for a semi-algebraic subset
$Y\subset \P J^k_{\C^n,0} \subset \P^* J^k_{\C^n,0}$.
\end{definition}

\begin{definition}
For a cylinder subset $X\subset \P{\cal O}_{\C^n,0}$ ($X=\pi_k^{-1}(Y)$,
$Y\subset \P J^k_{\C^n,0}$) its Euler characteristic $\chi(X)$ is defined
as the Euler characteristic $\chi(Y)$ of the set $Y$.
\end{definition}

Let $\psi: \P{\cal O}_{\C^n,0} \to A$ be a function with values in
an Abelian group $A$.

\begin{definition}
We say that the function $\psi$ is cylindric if, for each $a\ne 0$,
the set $\psi^{-1}(a)\subset \P{\cal O}_{\C^n,0}$ is cylindric.
\end{definition}

\begin{definition}
The integral of a cylindric function $\psi$ over the space
$\P{\cal O}_{\C^n,0}$ with respect to the Euler characteristic is
$$
\int_{\P{\cal O}_{\C^n,0}}\psi d\chi :=
\sum_{a\in A, a\ne 0} \chi(\psi^{-1}(a))\cdot a
$$
if this sum has sense in $A$. If the integral exists (has sense)
the function $\psi$ is said to be integrable.
\end{definition}

\begin{remarks}
{\bf 1.} In a similar way one can define the notion of the generalized
Euler characteristic $[X]$ of a cylindric subset of the space
$\P{\cal O}_{\C^n,0}$ (or of the space ${\cal O}_{\C^n,0}$)
with values in the Grothendieck ring of complex algebraic varieties
localized by the class $\L$ of the complex line and thus the
corresponding notion of integration (see, e.g., {\cite{C}}). For that one
should define $[X]$ as $[Y]\cdot\L^{-{n+k\choose k}}$.\newline
{\bf 2.} There are the same notions (of the Euler characteristic, of the
generalized Euler characteristic, and of the integral with respect to
the Euler characteristic) in the real setting as well. Since the Euler
characteristic $\chi(\L_\R)$ of the real line $\L_\R$ is equal to
$(-1)$, one has to define the Euler characteristic $\chi(X)$ of
a cylinder subset $X\subset \P{\cal E}_{\R^n,0}$ ($X=\pi_k^{-1}(Y)$,
$Y\subset \P J^k_{\R^n,0}$) as $(-1)^{-{n+k\choose k}}\cdot\chi(Y)$.
\end{remarks}

Let $\Z[[t]]$ (respectively $\Z[[t_1, \ldots, t_r]]$) be the group
(with respect to the addition) of formal power series in the variable
$t$ (respectively in the variables $t_1$, \dots, $t_r$).
As usual, for $\vv=(v_1, \ldots, v_r)\in \Z_{\ge 0}^r$,
$\tt^\vv=t_1^{v^1}\cdot\ldots\cdot t_r^{v_r}$; we assume $t^\infty=0$.

\begin{theorem}\label{theo5}
For each $\vv \in \Z_{\ge 0}^r$, the subset $\{g\in \P{\cal O}_{\C^2, 0}:
\vv(g)=\vv\}$ of $\P{\cal O}_{\C^2, 0}$ is cylindric.
Therefore the functions $\tt^{\vv(g)}$ and $t^{v(g)}$ $(v(g)=
\Vert \vv(g) \Vert = v_1(g)+\ldots+v_r(g))$ on $\P{\cal O}_{\C^2, 0}$ with
values in $\Z[[t_1, \ldots, t_r]]$ and $\Z[[t]]$ respectively are cylindric.
\end{theorem}

{\bf Proof} follows from the fact that, for $g\in {\goth m}^s$, $v_i(g)\ge s$,
i.e., the power series decomposition of $g\circ\varphi_i(\tau_i)$ starts
from terms of degree at least $s$. Therefore the functions $\tt^{\vv(g)}$
and $t^{v(g)}$ on $\P{\cal O}_{\C^2, 0}$ are integrable (since
$\sum\limits_{\vv\in\Z_{\ge 0}^r}\ell(\vv)\tt^\vv\in \Z[[t_1, \cdots, t_r]]$
for any integers $\ell(\vv)$). {$\Box$}

\begin{theorem} \label{theo6}
For $r>1$, $$\int\limits_{\P{\cal O}_{\C^2, 0}}\tt^{\vv(g)} d\chi =
\Delta^C(t_1,\ldots, t_r);$$
for $r\ge 1$, $$\int\limits_{\P{\cal O}_{\C^2, 0}}t^{v(g)} d\chi =
\zeta_C(t).$$
\end{theorem}

{\bf Proof} follows from Theorem~\ref{theo4}. For $\vv =(v_1, \ldots, v_r)
\in \Z_{\ge 0}^r$, let $k=1+\max\limits_{1\le i\le r}v_i$, and let
$Z_{\vv}\subset\P J^k_{\C^2, 0}$ be the set $\{j^kg\in\P J^k_{\C^2, 0}:
\vv(g)=\vv\}$ (we have mentioned that the condition $\vv(g)=\vv$ is
determined by the $k$--jet $j^kg$ of the germ $g$). The natural map
$Z_\vv\to \P F_\vv$ of the set $Z_\vv$ to the fibre $\P F_\vv$ of the
projectivization $\P{\hat S}_C$ of the extended semigroup ${\hat S}_C$ is
a locally trivial fibration the fibre of which is a complex vector space
of some dimension. Now Theorem~\ref{theo6} follows from the facts that
$\chi(\P{\hat S_C})=\Delta^C(t_1,\ldots, t_r)$,
$\zeta_C(t)=\Delta^C(t, \ldots, t)$. {$\Box$}

\section{Polynomial functions in two variables and the zeta-function
of the monodromy transformation at infinity}
\label{sec6}

There exists a global analogue of the discussed results for a plane
algebraic curve with one branch at infinity (\cite{CDG3}). Thus one can
say that this is an analogue of the statement about the coincidence of
the Poincar\'e series of the coordinate ring and the zeta-function of the
classical monodromy transformation of an irreducible curve singularity
(i.e., for $r=1$).

Let $F(x, y)$ be a complex polynomial of degree $d$. It is well
known that the map $F:\C^2\to\C$ is a locally trivial fibration
over the complement to a finite set in the target $\C$. Its fibre is
a generic level curve $C_{gen}$ of the polynomial $F$. Let
$h^\infty_F:C_{gen}\to C_{gen}$ be the monodromy transformation of this
fibration corresponding to the loop $\gamma(\tau)=R\cdot exp(2\pi i\tau)$
with real $R$ large enough.
The transformation $h^\infty_F$ is called the monodromy transformation
of the polynomial $F$ at infinity. Let $\zeta^\infty_F(t)$ be the
zeta--function of the monodromy transformation $h^\infty_F$.

Let $C_\lambda={(x, y)\in\C^2: F(x, y)=\lambda}$ be the level curve of
the polynomial $F$.
Assume that the curve $C_0$ has only one branch at infinity, i.e., that
the the projective curve $\overline{C_0}$ has one point $A_{\infty}$
at infinity and that it has only one branch at this point. This implies
that the curve $\overline{C_{\lambda}}$ has only one branch at infinity
(i.e., at the point $A_{\infty}$) for any $\lambda\in \C$.

Without loss of generality one can assume that $A_{\infty} = (0:1:0)$. The
germ of the curve $\overline{C_0}$ at the point $A_{\infty}$ is determined
by the local equation $\{Q(u,v)=0\}$, where $Q(u,v) = u^d\cdot P(v/u,1/u)$;
the infinite line $L_{\infty}$ is given by the equation $\{u=0\}$.

To the curve $C_0$ there corresponds a semigroup $\Gamma\in\Z_{\ge0}$~---
the semigroup of (orders of) poles along the curve $\overline{C_0}$ of polynomials
in the affine plane $\C^2$ at the point $A_{\infty}$, i.e.,
$$
\Gamma := \{ - v(H(v/u,1/u)) | H(x,y)\in \C[x,y],
\mbox{ $H$ is not divisible by $F$}
\}\subset\Z.
$$

One can easily see that in fact $\Gamma\subset\Z_{\ge0}$, since a
polynomial not identically equal to zero on the curve $C_0$ cannot
have a zero at the point $A_{\infty}$ of it.
Let $P_\Gamma(t)$ be the Poincar\'e series for the semigroup $\Gamma$:
$$P_\Gamma(t)=\sum\limits_{i\in\Gamma}t^i.$$

\begin{theorem}\label{theo7}
$P_\Gamma(t)= \zeta^{\infty}_F(t)$.
\end{theorem}

The {\bf proof} is essentially the same as of Theorem~\ref{theo1}. Using
known combinatorial properties of the semigroup $\Gamma$ (given
by the Abhyankar--Moh theorem: \cite{AM1}, \cite{AM2}) one can
compute the Poincar\'e series $P_\Gamma(t)$ in terms of the minimal
embedded resolution of the plane curve $C_0$.
The zeta-function $\zeta^{\infty}_F(t)$ of the polynomial is given
by an analogue of the A'Campo formula. The results of the computations
show that they coincide. {$\Box$}

\begin{remark}
One can easily see that the statement of Theorem~\ref{theo7} can be
rewritten in the form
$$
\zeta^{\infty}_F(t)=\int\limits_{\P\C[x, y]}t^{-v(H)}d\chi,
$$
where $v(H)=v(H(v/u, 1/u))$. The integral with respect to the Euler
characteristic over the projectivization $\P\C[x, y]$ of the space
$\C[x, y]$ of polynomials in two variables is defined in the same way
as one in the section~\ref{sec5}.
\end{remark}

\bigskip
\noindent Translated by S.M.Gusein-Zade.


\begin{thebibliography}{12}
\bibitem[1]{A'C} A'Campo N. La fonction z\^eta d'une monodromie.
Comment. Math. Helv., 50, 233--248 (1975).
\bibitem[2]{AM1} Abhyankar S.S., T.T. Moh T.T. Newton-Puiseux
expansion and generalized Tschirnhaussen transformation I.
Cr\'elle, 260, 47-83 (1973).
\bibitem[3]{AM2} S.S. Abhyankar, T.T. Moh: Newton-Puiseux
expansion and generalized Tschirnhaussen transformation II.
Cr\'elle, 261, 29-54 (1973).
\bibitem[4]{C} Craw A. An introduction to motivic integration.
Preprint math.AG/9911179.
\bibitem[5]{DL} Denef J., Loeser F. Lefschetz numbers of iterates
of the monodromy and truncated arcs. Preprint math.AG/0001105.
\bibitem[6]{EN} Eisenbud D., Neumann W. Three-dimensional link theory and
invariants of plane curve singularities. Ann. of Math. Studies 110,
Princeton Univ. Press, Princeton, NJ, 1985.
\bibitem[7]{CDG1} Campillo A., Delgado F., Gusein-Zade S.M.
The extended semigroup of a plane curve singularity. Proceedings of
the Steklov Institute of Mathematics, v.221, 139--156 (1998).
\bibitem[8]{CDG2} Campillo A., Delgado F., Gusein-Zade S.M.
On the monodromy of a plane curve singularity and the Poincar\'e series
of its ring of functions. Functional Analysis and its Applications,
v.33, N 1, 66-68 (1999).
\bibitem[9]{CDG3} Campillo A., Delgado F., Gusein-Zade S.M.
On the monodromy at infinity of a plane curve and the Poincar\'e series
of its coordinate ring. In: Proceedings of the International Conference
devoted to the 90-th birthday of L.S.Pontryagin (Moscow, 31.8-06.9 1998),
v.7 (Geometry and Topology), p.49--54; VINITI: Itogi nauki i tekhniki
(Results of science and techniques), Series "Contemporary mathematics
and its applications", v.68, Moscow, 1999 (in Russian).
\bibitem[10]{CDG4} Campillo A., Delgado F., Gusein-Zade S.M. The Alexander
polynomial of a plane curve singularity, and the ring of functions on the
curve. Russian Math. Surveys, 54 (1999), no. 3 (327), 157--158.
\bibitem[11]{ZT} Zariski, O. \& Teissier, B. Le probl\`eme des
modules pour les branches planes. Hermann, Paris 1986
\end{thebibliography}
\end{document}